\documentclass[11pt,a4paper]{article}

\usepackage{amsmath}     
\usepackage{dsfont}		
\usepackage{amsfonts}		
\usepackage{amsthm}
\usepackage{graphicx}
\usepackage{amssymb}

\usepackage{latexsym}
\usepackage{euscript,makeidx,color,mathrsfs,algorithm}
\usepackage{enumerate}
\usepackage{tabu}

\usepackage[colorlinks,linkcolor=blue,anchorcolor=green,citecolor=red]{hyperref}

\def\sqr#1#2{{\vcenter{\vbox{\hrule height.#2pt
              \hbox{\vrule width.#2pt height#1pt \kern#1pt \vrule width.#2pt}
              \hrule height.#2pt}}}}
\def\5n{\negthinspace \negthinspace \negthinspace \negthinspace \negthinspace }
\def\4n{\negthinspace \negthinspace \negthinspace \negthinspace }
\def\3n{\negthinspace \negthinspace \negthinspace }
\def\2n{\negthinspace \negthinspace }
\def\1n{\negthinspace }

\def\dbF{\mathbb{F}}

\def\dbK{\mathbb{K}}

\def\dbN{\mathbb{N}}

\def\dbP{\mathbb{P}}

\def\dbR{\mathbb{R}}

\def\dbU{\mathbb{U}}

\def\dbX{\mathbb{X}}

\def\={\buildrel \triangle \over =}

\def\ds{\displaystyle}

\def\d{\delta}

\def\l{\lambda}

\def\t{\tau}
\def\f{\varphi}

\def\D{\Delta}

\def\L{\Lambda}
\def\Si{\Sigma}

\def\O{\Omega}

\def\me{\mathbb{E}}
\def\rd{\,\mathrm d}  
\def\bal{\begin{aligned}}
\def\eal{\end{aligned}}

\def\cC{{\cal C}}

\def\cF{{\cal F}}

\def\cJ{{\cal J}}

\def\no{\noindent}

\def\ms{\medskip}

\def\q{\quad}
\def\qq{\qquad}

\def\lt{\left}
\def\rt{\right}
\def\lan{\langle}
\def\ran{\rangle}
\def\llan{\left\langle}
\def\rran{\right\rangle}

\def\rf{\eqref}

\def\h{\widehat}
\def\wt{\widetilde}

\def\cd{\cdot}
\def\cds{\cdots}

\def\ae{\hbox{\rm a.e.}}

\def\deq{\triangleq}

\def\({\Big (}
\def\){\Big )}
\def\[{\Big[}
\def\]{\Big]}

\def\bde{\begin{definition}\label}
\def\ede{\end{definition}}
\def\be{\begin{equation}}
\def\bel{\begin{equation}\label}
\def\ee{\end{equation}}
\def\beq{\begin{equation*}\begin{aligned}}
\def\eeq{\end{aligned}\end{equation*}}
\def\bt{\begin{theorem}\label}
\def\et{\end{theorem}}
\def\bc{\begin{corollary}\label}
\def\ec{\end{corollary}}
\def\bl{\begin{lemma}\label}
\def\el{\end{lemma}}
\def\bp{\begin{proposition}\label}
\def\ep{\end{proposition}}
\def\bas{\begin{assumption}\label}
\def\eas{\end{assumption}}
\def\br{\begin{remark}\label}
\def\er{\end{remark}}
\def\bex{\begin{example}\label}
\def\ex{\end{example}}
\def\ba{\begin{array}}
\def\ea{\end{array}}
\def\ed{\end{document}}

\def\square#1{\vbox{\hrule\hbox{\vrule height#1%
     \kern#1\vrule}\hrule}}
\def\rectangle#1#2{\vbox{\hrule\hbox{\vrule height#1%
     \kern#2\vrule}\hrule}}

\font\tenbb=msbm10 \font\sevenbb=msbm7 \font\fivebb=msbm5

\newfam\bbfam
\scriptscriptfont\bbfam=\fivebb \textfont\bbfam=\tenbb
\scriptfont\bbfam=\sevenbb

\newtheorem{theorem}{\hskip 1.3em Theorem}[section]
\newtheorem{definition}[theorem]{\hskip 1.3em Definition}
\newtheorem{proposition}[theorem]{\hskip 1.3em Proposition}
\newtheorem{corollary}[theorem]{\hskip 1.3em Corollary}
\newtheorem{lemma}[theorem]{\hskip 1.3em Lemma}
\newtheorem{remark}[theorem]{\hskip 1.3em Remark}
\newtheorem{example}[theorem]{\hskip 1.3em Example}

\newtheorem{assumption}[theorem]{\hskip 1.3em Assumption}

\allowdisplaybreaks 

\makeatletter
   
   \@addtoreset{equation}{section}
\makeatother

\begin{document}

\title{Convergence rates of a discrete feedback control arising in mean-field linear quadratic 
optimal control problems\thanks{This work is supported by
the National Natural Science Foundation of China (11801467).}}

\author{Yanqing Wang\thanks{
School of Mathematics and Statistics, Southwest University, Chongqing 400715, China.  {\small\it
e-mail:} {\small\tt yqwang@amss.ac.cn}. \ms}}

\date{Feb. 3, 2023}
\maketitle

\begin{abstract}
In this work, we propose a feedback control based temporal discretization for linear 
quadratic optimal control problems (LQ problems)
governed by controlled mean-field stochastic differential equations. 
We firstly decompose the original problem into two problems: a stochastic LQ problem and a
deterministic one.
Secondly, we discretize both LQ problems one after another relying on Riccati equations and control's feedback representations.
 Then, we prove the convergence rates for the proposed discretization and present an effective algorithm.
Finally, a numerical example is provided to support the theoretical finding.
\end{abstract}

\ms

\no\bf Keywords: \rm Convergence rate, mean-field linear quadratic optimal control problem,  closed--loop optimal control strategy,
 Riccati equation


\ms

\no\bf AMS 2010 subject classification: \rm 
 93E20,
 49N10,
 93B52

\section{Introduction}

The present work is concerned with the numerical scheme for the linear quadratic optimal control problems 
(LQ problems, for short) governed
by controlled mean-field stochastic differential equations (MF-SDEs, for short), which are special stochastic optimal 
control problems. 

In 1938, the idea of mean field interaction in mathematical physics was proposed 
by Vlasov. From 1950s, a kind of MF-SDEs, named by McKean-Vlasov SDEs, was proposed and then systematically studied. Up to now, theory of MF-SDEs has significant impact in several areas such as physics, financial mathematics, stochastic optimal control, social interactions.

As three well-recognized fundamental results of stochastic optimal control theory, Pontryagin-type maximum principle 
(MP, for short) (see, e.g., \cite{Boltyanskii-Gamkrelidze-Pontryagin56, Peng90}), Bellman dynamic programming principle (DPP, for short) (see, e.g., \cite{Bellman57, Crandall-Ishii-Lions92,Pham16}) and LQ theory (see, e.g., \cite{Kalman60,Kushner65}) have been established in the past
half century.
For the stochastic optimal control problems subject to MF-SDEs, there still exists some works; see, e.g., 
\cite{Ahmed07,Buckdahn-Djehiche-Li11,Andersson-Djehiche11,MeyerBrandis-Oksendal-Zhou12,Yong13SICON,Huang-Li-Yong15,Li-Sun-Xiong19}.


However, general stochastic optimal control problems, except for some special  case 
such as one-dimensional state equation,
can not admit explicitly closed form solutions. Therefore, effective numerical schemes for these problems are necessary,
and attract some attentions. 
Up to now, there are two widely applied and popular approaches to numerically solve stochastic optimal control problems:
DPP based approach and MP based one.
The first approach relies on DPP and Hamilton-Jacobi-Bellman 
(HJB, for short) equation. By establishing the bridge between the original problem and an HJB equation,
then numerically solving the HJB equation, one can obtain an approximation of the optimal control (see, e.g.,
\cite{Feng-Glowinski-Neilan13,Feng-Jensen17,Beck-E-Jentzen19}). However, due to the nonlinearity of the HJB equation, 
this approach is computational expensive.
For the second approach,
thanks to MP, the optimal pair satisfies a coupled 
forward-backward stochastic differential equation (FBSDE, for short).
The difficulty and complexity of the second approach stem from the fact that this equation is coupled and  the solution to
this equation should be adapted. To decouple this FBSDE,
 a (stochastic) gradient descent algorithm can be adopted (see, e.g., \cite{Dunst-Prohl16,Archibald-Bao-Yong20,Archibald-Bao-Yong-Zhou20,Prohl-Wang21a,Prohl-Wang22}). Besides,
gradient projection methods for stochastic optimal control problems with deterministic control variables are also
proposed; see \cite{Du-Shi-Liu13,Gong-Liu-Tang-Zhao-Zhou17}.
To guarantee the solution's adaptivity, conditional expectations are introduced (see, e.g., \cite{ZhangJF04,Bender-Denk07}). Although there exists numerous works to simulate conditional expectations (see, e.g., \cite{Gobet-Lemor-Warin05,Bender-Denk07,Hu-Nualart-Song11,Wang-Zhang11,Briand-Labart14}), the study of computation on conditional expectations is
far from mature till now.
From this perspective, the computational cost of  MP based approach is high.

To reduce the computation cost, recently based on LQ theory, Riccati based schemes for LQ problems subject to 
controlled SDEs are proposed in \cite{Wang21Riccati}. Avoiding solving coupled equations and 
simulating conditional expectations,
these schemes can greatly save computation time.

\ms

In this work, adopting LQ theory of 
MF-SDEs and Riccati equations, we aim at proposing an efficient  algorithm with convergence rates to 
solve LQ problems governed by MF-SDEs, which will be called MF-LQ problems. To be specific, 
Let $(\Omega, {\mathcal F}, {\mathbb F}, {\mathbb P})$ be a complete filtered probability space on which a 
standard one-dimensional Brownian motion $W=\{W(t):\,t\geq 0\}$ is defined, where $\dbF=\{\cF_t\}_{t\geq 0}$ is the natural filtration of $W$ augmented by all the $\dbP$-null sets in $\cF$.
Consider the following cost functional 
\beq
\cJ( u(\cd))&=\me\Big[ \int_0^T\big( \lan Qx(t), x(t)\ran+\lan \bar Q\me[x(t)], \me[x(t)]\ran+ \lan Ru(t), u(t)\ran+\lan \bar R\me[u(t)], \me[u(t)]\ran \big)\rd t\\
&\q
+\lan G x(T), x(T)\ran +\lan \bar G \me[x(T)],\me[x(T)]\ran \Big]
\eeq
subject to a controlled MF-SDE:
\bel{sde1}
\lt\{
\bal
&\rd x(t)=\lt(A x(t)+\bar A \me[x(t)]+Bu(t)+\bar B \me[u(t)]\rt)\rd t\\
&\qq\qq+\lt(C x(t)+\bar C \me[x(t)]+D u(t)+\bar D\me[u(t)]\rt)\rd W(t) \q &\forall\, t \in [0,T]\,,\\
&x(0)=x_0\in \dbR^n\,,
\eal
\rt.
\ee
where $A,\,\bar A,\,C,\,\bar C,\,Q,\,\bar Q,\,G,\,\bar G\in\dbR^{n\times n}$, $B,\,\bar B,\,D,\,\bar D\in\dbR^{n\times m}$, $R,\,\bar R\in\dbR^{m\times m}$, and $u(\cd)$ is the control variable. The space of admissible controls for \rf{sde1} is
\beq
L^2_\dbF(0,T;\dbR^m)\deq\Big\{
\f :(0,T) \times \Omega \to\dbR^m \,\big|\,
\f(\cd)\hbox{ is $\dbF$-adapted and } \int_0^T\me[\|\f(t)\|^2]\rd t<\infty\Big\}\,.
\eeq
The considered MF-LQ problem can be stated as follows.\\
\no {\bf Problem (MF-LQ)}. For given $x_0\in\dbR^n$, find an optimal control $u^*(\cd)\in L^2_\dbF(0,T;\dbR^m)$ such that
\beq
\cJ(u^*(\cd))=\inf_{u(\cd)\in L^2_\dbF(0,T;\dbR^m)}\cJ(u(\cd))\,.
\eeq
Our goal is to numerically approximate the optimal stochastic control 
$u^*(\cd)$ of Problem (MF-LQ).

%
 
The main strategy for discretization of Problem (MF-LQ)  is as 
follows: Firstly, thanks to the optimal control's state feedback form and two corresponding differential Riccati equations,  
we introduce a stochastic LQ problem, Problem (LQ)$^{1,t}$, related to the first Riccati equation, and 
a discrete stochastic LQ problem, Problem (LQ)$^{1,t}_\t$, to approximate Problem (LQ)$^{1,t}$;
Secondly, relying on the solution to the first aforementioned differential Riccati equation, we introduce a deterministic LQ problem,
Problem (LQ)$^{2,t}$, and its temporal discretization, Problem (LQ)$_\t^{2,t}$, whose related difference Riccati equation is proved
to converge to the solution of the second aforementioned differential Riccati equation;
Thirdly, by inverting the above two steps and adopting the optimal controls' state feedback representations of the above
four auxiliary  LQ problems, we  can propose a discretization for Problem (MF-LQ), prove the convergence of our proposed discretization and derive the rates of convergence. This is one contribution of this work.

As an application of our scheme, we actually can propose efficient and accurate algorithms with rates to numerically solve differential Riccati equations
and a family of coupled mean-field forward-backward stochastic differential equations.
This is another contribution of this work.



The rest of this paper is organized as follows. In Section \ref{pre}, we introduce some notations, review
results on solvability of Problem (MF-LQ), present temporal discretization of stochastic LQ problems and their corresponding  
Riccati equations, and state the convergence with rates for aforementioned discretization of Riccati equation which is
related to the optimal control's state feedback representation. 
In Section \ref{rate}, 
by virtue of closed--loop optimal control strategy, we 
prove the convergence rates of temporal discretization for Problem (MF-LQ),  and propose an accurate and efficient algorithm, as well as carry out an 
example to verify our theoretical results.
%

\section{Preliminaries}\label{pre}

\subsection{Notations and solvability of Problem (MF-LQ)}\label{not1}

Let $\dbK$ be $\dbR^n$ or $\dbR^m$. 
Denote 
 by $L^2_\dbF(0,T;\dbK)$ the space of $\dbK$-valued, $\dbF$-adapted processes $\f(\cd)$ satisfying $\me\big[\int_0^T \|\f(t)\|_\dbK^2\rd t\big]<\infty$; 
 by $L^2_\dbF(\O;C([0,T];\dbK))$ the 
subspace of $L^2_\dbF(0,T;\dbK)$, such that every element $\f(\cd)$ has a continuous path and satisfies 
$\me\big[\sup_{t\in[0,T]}\|\f(t)\|_\dbK^2 \big]<\infty$;
by $L^2(0,T;\dbK)$ the space of $\dbK$-valued functions $\f(\cd)$ satisfying $\int_0^T \|\f(t)\|_\dbK^2\rd t<\infty$.
For any matrix $M$, denote $\|M\| \deq \sqrt{\rho(M^\top M)}$, where $\rho$ represents the spectral radius. 
A symmetric matrix $M>0(\geq 0)$ means that it is positive 
definite (positive semi-definite).
Denote by $I_n$ the identity matrix of size $n$.

We denote by $I_\tau = \{ t_{n}\}_{n=0}^N \subset [0,T]$ a time mesh with maximum step size $\tau :=\max\{t_{n+1}-t_n:\,n=0,1,\cds,N-1\}$, and
$\D_nW=W(t_n)-W(t_{n-1})$ for all $n=1,\cds,N $. 
For a given time mesh $I_\t$, 
we can define $\nu: [0,T]\to I_\t $ 
by
\bel{w128e1}
\bal
&\nu (t)=t_k\qq \forall\, t\in [t_k,t_{k+1}),\,k=0,1,\cds,N-1\,.
\eal
\ee
For simplicity, we choose a uniform partition, i.e.  step-size $\t=T/N$ and $\t\leq 1$. 
Denote by $\cC$ a generic positive constant
independent of data and $I_\t$,
which may vary from one place
to another.


\ms

Throughout this paper, to guarantee the solvability of Problem (MF-LQ),  we make the following assumption.\\
{\bf (A)} $Q\geq 0,\,Q+\bar Q\geq0,\,G\geq0,\,G+\bar G\geq 0,\,R> 0,\,R+\bar R>0$.

\ms

For simplicity,
denote $\h \varphi \deq \varphi+\bar \varphi$, for $\varphi=A,\,B,\,C,\,D,\,Q,\,R,\,G$.

The results in this work can be extended to 
cases with a higher-dimensional $W(\cd)$,  a general partition $I_\t$, and/or time-varying deterministic coefficients without significant difficulties.
Note that for the time-varying coefficients, the assumption 
should be changed to the following:\\
{\bf (A1)} $Q(\cd)\geq 0,\,\h Q(\cd)\geq0,\,G\geq0,\,\h G\geq 0,\,R(\cd)-\d I_m\geq 0,\,\h R(\cd)-\d I_m\geq0$, where $\d$ 
is a positive constant;\\
{\bf (A2)} $\varphi(\cd),\,\bar \varphi(\cd)$ are Lipchitz continuous, for $\varphi(\cd)=A(\cd),\,B(\cd),\,C(\cd),\,D(\cd),\,Q(\cd),\,R(\cd)$.

\ms
The following result is on the regularity of the state $x(\cd)$; see \cite[Proposition 2.6]{Yong13SICON}.

\bl{w1003l1}
Let {\rm(A)} hold. Then for any $(x_0,u(\cd))\in\dbR^n\times L^2_\dbF(0,T;\dbR^m)$, equation \rf{sde1} admits a unique
solution $x(\cd)\in L^2_\dbF(\O;C([0,T];\dbR^n))$.

\el

The following result states Pontryagin-type maximum principles for Problem (MF-LQ),
which guarantees solvability of Problem (MF-LQ). The reader can refer to  \cite[Proposition 2.8, Theorems 3.1 \& 3.2]{Yong13SICON}.

\bl{MP-1}
Under the assumption {\rm (A)},
Problem {\rm(MF-LQ)} admits a unique optimal pair
$(x^*(\cd), u^*(\cd))\in L^2_\dbF(0,T;\dbR^n) \times L_\dbF^2(0,T;\dbR^m)$, and the unique solvability of Problem 
{\rm(MF-LQ)} is
equivalent to the unique solvability of the following equation:
\bel{w1109e1}
\lt\{
\bal
& \rd x^*(t)=\lt(A x^*(t)+\bar A \me[x^*(t)]+Bu^*(t)+\bar B \me[u^*(t)]\rt)\rd t\\
&\qq\qq+\lt(C x^*(t)+\bar C \me[x^*(t)]+D u^*(t)+\bar D\me[u^*(t)]\rt)\rd W(t) \qq \forall\, t \in [0,T]\,,\\
&\rd y(t)=-\big(A^\top y(t)+\bar A^\top \me[y(t)]+ C^\top z(t)+\bar C^\top \me[z(t)]- Qx^*(t)-\bar Q \me[x^*(t)] \big)\rd t\\
&\qq\qq+z(t)\rd W(t)  \qq \forall\, t\in [0,T] \,,\\
&x^*(0)=x_0\, ,\q y(T)=-Gx^*(T)\,,
\eal
\rt.
\ee
 together with the optimality condition
\beq
R u^*(\cd)+\bar R \me[u^*(\cd)]-\big(B^\top y(\cd)+\bar B^\top \me[y(\cd)]+D^\top z(\cd)+\bar D^\top \me[z(\cd)]\big)=0 \q \ae 
\eeq
\el

In the following, we tend to state the closed--loop optimal control strategy of MF-LQ problems, where two differential Riccati equations 
are needed.
By \cite{Yong13SICON},
these two Riccati equations (with time variable suppressed) related to Problem (MF-LQ) read
\bel{Riccati-1}
\lt\{
\bal
&P'+PA+A^\top P+C^\top PC+Q\\
&\qq-(B^\top P+D^\top PC)^\top \Si_0^{-1}( B^\top P+D^\top PC)=0\q \mbox{in } [0,T]\,,\\
&P(T)=G\,,\\
\eal
\rt.
\ee
and
\bel{Riccati-2}
\lt\{
\bal
&\Pi'+\Pi \h A+\h A^\top \Pi+\h C^\top P\h C+\h Q\\
&\qq-(\h B^\top \Pi+\h D^\top P\h C)^\top \Si_1^{-1}
(\h B^\top \Pi+\h D^\top P\h C)=0\q \mbox{in } [0,T]\,,\\
&\Pi(T)=\h G\,,\\
\eal
\rt.
\ee
where 
\bel{w929e1}
\Si_0=R+D^\top PD\,,\qq\Si_1=\h R+\h D^\top P\h D\,.
\ee
The following result on solvability of these two equations and the optimal control's closed--loop representation
 belongs to
 \cite[Theorem 4.1]{Yong13SICON}.

\bl{Yong-Riccati}
Under the assumption {\rm (A)},
Riccati equations \rf{Riccati-1} and \rf{Riccati-2} admit unique solutions $P(\cd)$ and $\Pi(\cd)$ which are positive
semi-definite. Further, the optimal
control to Problem \rm{(MF-LQ)} has the following feedback representation (with time variable suppressed):
\bel{w916e1}
u^*=-\Si_0^{-1}\big[B^\top P+D^\top P C\big](x^*-\me[x^*])-\Si_1^{-1}\big[\h B^\top\Pi+\h D^\top P\h C \big]\me[x^*]\,.
\ee
%
\el
%


\subsection{Temporal discretization of parametrized MF-LQ problems} \label{error-of-LQ}

We employ a uniform mesh $I_\tau=\{t_i\}_{i=0}^N$ covering $[0,T]$, and for any $t_l\in I_\t,\,l=0,1,\cds,N-1$, consider step size process pair
$(x_\t(\cd), u_\t(\cd)) \in \dbX_{\t}(t_l,T) \times \dbU_{\t}(t_l,T)\subset L^2_\dbF(t_l,T;\dbR^n )\times L_\dbF^2(t_l,T;\dbR^m)$, 
where
\beq
\dbX_{\t}(t_l,T) &\deq \lt\{x(\cd) \in L^2_\dbF(t_l,T;\dbR^n )\,|\, x(t)=x(t_n), \,\,   \forall t\in [t_n, t_{n+1}),  \,\, n=l,l+1,\cds, \, N-1\rt\}\, ,\\
\dbU_{\t}(t_l,T) &\deq \lt\{u(\cd) \in L^2_\dbF(t_l,T;\dbR^m)\,|\, u(t)=u(t_n),  \,\,    \forall t\in [t_n, t_{n+1}),  \,\, n=l,l+1,\cds, \, N-1\rt\}\, ,
\eeq
and for any $x(\cd)\in \dbX_{\t}(t_l,T)$  and $u(\cd)\in \dbU_{\t}(t_l,T)$,
\beq
\|x(\cd)\|_{\dbX_{\t}(t_l,T)}\deq\Big(\t \sum_{n=l}^{N-1}\me\big[\|x(t_n)\|^2\big]\Big)^{1/2}\,,\quad
\mbox{and} \quad \|u(\cd)\|_{\dbU_{\t}(t_l,T)}\deq\Big(\t \sum_{n=l}^{N-1}\me\big[\|u(t_n)\|^2\big]\Big)^{1/2}.
\eeq
We also need the following two deterministic spaces:
\beq
\bar\dbX_{\t}(t_l,T) &\deq \lt\{x(\cd) \in L^2(t_l,T;\dbR^n )\,| \, x(t)=x(t_n), \,\,   \forall t\in [t_n, t_{n+1}),  \,\, n=l,l+1,\cds, \, N-1\rt\}\, ,\\
\bar\dbU_{\t}(t_l,T) &\deq \lt\{u(\cd) \in L^2(t_l,T;\dbR^m)\,|\, u(t)=u(t_n),  \,\,    \forall t\in [t_n, t_{n+1}),  \,\, n=l,l+1,\cds, \, N-1\rt\}\, ,
\eeq
and for any $x(\cd)\in \bar\dbX_{\t}(t_l,T)$  and $u(\cd)\in \bar\dbU_{\t}(t_l,T)$,
\beq
\|x(\cd)\|_{\bar\dbX_{\t}(t_l,T)}\deq\Big(\t \sum_{n=l}^{N-1} \|x(t_n)\|^2\Big)^{1/2}\,,\quad
\mbox{and} \quad \|u(\cd)\|_{\bar\dbU_{\t}(t_l,T)}\deq\Big(\t \sum_{n=l}^{N-1} \|u(t_n)\|^2 \Big)^{1/2}.
\eeq

For the optimal pair $(x^*(\cd),u^*(\cd))$ of Problem (MF-LQ), by Lemma \ref{w1003l1}, it follows that 
$x^*(\nu(\cd))\in \dbX_\t(0,T)$, where $\nu(\cd)$ is defined in \rf{w128e1}. Then by Lemma \ref{Yong-Riccati}, we know that $u^*(\nu(\cd))\in\dbU_\t(0,T)$.
By H\"older's inequality, it is easy to see that for any $(x(\cd),u(\cd))\in \dbX_\t (t_l,T)\times \dbU_\t (t_l,T), l=0,1,\cds,N-1$,
$(\me[x(\cd)],\me[u(\cd)])\in \bar \dbX_\t (t_l,T)\times \bar\dbU_\t (t_l,T)$.


Now, we introduce parametrized MF-LQ problems.  For any $t\in[0,T]$, $l=0,1,\cds,N-1$, consider the following systems:
\bel{mf-sde}
\lt\{
\bal
&\rd x(s)=\lt(A x(s)+\bar A \me[x(s)]+Bu(s)+\bar B \me[u(s)]\rt)\rd s\\
&\qq\qq+\lt(C x(s)+\bar C \me[x(s)]+D u(s)+\bar D\me[u(s)]\rt)\rd W(s) \q &\forall\, s \in [t,T]\,,\\
&x(t)=x_t\in \dbR^n\,,
\eal
\rt.
\ee
and
\bel{mf-sdde}
\lt\{
\bal
& x_{\t}(t_{k+1})-x_{\t}(t_k)= \tau \big[A x_{\t}(t_{k})+\bar A \me[x_\t(t_k)]
 +Bu_{\t}(t_k)+\bar B \me[u_\t(t_k)] \big]\\
 &\qq\qq\qq\qq\qq+ \big[Cx_{\t}(t_k)+\bar C\me[x_\t(t_k)]+D u_\t (t_k)+\bar D\me[u_\t(t_k)]\big] \D_{k+1}W\\
& \qq\qq\qq\qq\qq\qq\qq\qq\qq\qq \q k=l,l+1,\cds,N-1\, ,\\
& x_{\t}(t_l)= x_{t_l}\in\dbR^n\,,
\eal
\rt.
\ee
as well as the cost functionals
\bel{cost-t}
\bal
&\cJ(t,x_t; u(\cd))\\
&=\me\Big[ \int_t^T\big( \lan Qx(s), x(s)\ran+\lan \bar Q\me[x(s)], \me[x(s)]\ran+ \lan Ru(s), u(s)\ran+\lan \bar R\me[u(s)], \me[u(s)]\ran \big)\rd s\\
&\q
+\lan G x(T), x(T)\ran +\lan \bar G \me[x(T)],\me[x(T)]\ran \Big]
\eal
\ee
and
\bel{cost-l}
\bal
\cJ_\t(t_l,x_{t_l}; u_\t(\cd))
&=  \|\sqrt{Q}x_\t(\cd)\|^2_{\dbX_\t(t_l,T)} +\lan \bar Q \me[x_\t(\cd)], \me[x_\t(\cd)]\ran_{\bar\dbX_\t (t_l,T)}\\
&\q + \| \sqrt{R}u_\t(\cd)\|^2_{\dbU_\t(t_l,T)}+\lan \bar R \me[u_\t(\cd)], \me[u_\t(\cd)]\ran_{\bar\dbU_\t (t_l,T)}\\
&\q+ \me\lan G x_\t(T), x_\t(T)\ran +\lan \bar G \me[x_\t(T)], \me[x_\t(T)]\ran\,.
\eal
\ee
Here, $x_\t(\cd)$ and $u_\t(\cd)$ are defined as follows:
\beq
x_\t(t)=x_\t(t_k)\,,\q u_\t(t)=u_\t(t_k)\qq \forall\,t\in [t_k,t_k+1),\,k=l,l+1,\cds,N-1.
\eeq


In the below, we construct two auxiliary LQ problems.\\
\no {\bf Problem (LQ)$^{1,t}$}. For given $x_t\in\dbR^n$,  minimize $\cJ(t,x_t;u(\cd))$ over $L^2_\dbF(t,T;\dbR^m)$ subject to \rf{mf-sde} with $\bar A=\bar C=\bar G=0,\,\bar B=\bar D=0,\,\bar R=0$.\\
\no {\bf Problem (LQ)$_\t^{1,t_l}$}. For given $x_{t_l}\in\dbR^n$,  minimize $\cJ_\t(t_l,x_{t_l};u_\t(\cd))$ over $\dbU_\t (t_l,T)$ subject to \rf{mf-sdde} with $\bar A=\bar C=\bar G=0,\,\bar B=\bar D=0,\,\bar R=0$.

The following result comes from \cite[Theorem 4.3]{Zhou02}.
\bl{Zhou-Riccati}
Under the assumption {\rm (A)},
for any $l=0,1,\cds,N-1$,
Problem \rm{(LQ)$_\t^{1,t_l}$} admits a unique optimal control $u_\t^*(\cd)$, which has the following feedback representation:
\bel{w929e3}
u_\t^*(t_k)=-(W_k^1)^{-1}H_k^1x_\t^*(t_k) \qq \forall\,k=l,l+1,\cds,N-1\,,
\ee
and
\bel{w929e5}
\inf_{u_\t(\cd)\in \dbU_\t (t_l,T)}\cJ_\t(t_l,x_{t_l};u_\t(\cd))=\lan P_l x_{t_l},x_{t_l}\ran\,,
\ee
where
\bel{w929e4}
\bal
W_k^1=R+B^\top P_{k+1}B\t+D^\top P_{k+1}D\,,\qq
H_k^1=B^\top P_{k+1}(I_n+A\t)+D^\top P_{k+1}C\,,
\eal
\ee
and $P_\cd$ solve the following difference Riccati equations: 
\bel{Riccati-D1}
\lt\{
\bal
&P_k=\t Q+(I_n+A\t)^\top P_{k+1}(I_n+A\t)+\t C^\top P_{k+1}C-\t(H_k^1)^\top (W_k^1)^{-1}H_k^1\\
&\qq\qq\qq\qq\q k=0,1,\cds,N-1,\\
&P_N=G\,.
\eal
\rt.
\ee
\el

Problem (LQ)$^{1,t}$ and Problem (LQ)$_\t^{1,t_l}$ are standard stochastic LQ problems, and under the assumption 
(A), they are uniquely solvable; see \cite{Yong-Zhou99,Zhou02}. Besides, based on the LQ theory, the following result holds
(see  \cite{Wang21Riccati}). 

\bl{w229l4}
Let $P(\cd)\,,P_\cd$ be the solutions to Riccati equations \eqref{Riccati-1}, \rf{Riccati-D1} respectively.
Then there exists a positive constant $\cC$ 
independent of $\t$ such that
\bel{w301e7}
\lt\{
\bal
& \sup_{t\in[0,T]}\|P(t)\| +\max_{l=0,1,\cds,N}\|P_l\| \leq \cC\,, \\
& \|P(t)-P(s)\| \leq \cC|t-s|\qq  \forall\,t,s\in [0,T]\,,\\
&\max_{k=0,1,\cds,N}\|P(t_k)-P_k\|\leq \cC \sqrt{\t}\,.
\eal
\rt.
\ee
\el

Note that Lemma \ref{w229l4} is crucial to derive the convergence rate of the discretization for Problem (LQ)$^{1,t}$.
Actually, we can conclude that
\beq
\sup_{k=0,1,\cds,N}\me\big[\|x^*(t_k)-x^*_\t(t_k)\|^2+\|u^*(t_k)-u^*_\t(t_k)\|^2 \big] \leq \cC \t\|x_0\|^2\,,
\eeq
where $(x^*(\cd),u^*(\cd))$ and $(x^*_\t(\cd),u^*_\t(\cd))$ are the optimal pairs of Problem (LQ)$^{1,0}$ and Problem
(LQ)$_\t^{1,0}$ respectively; see \cite[Theorem 3.2]{Wang21Riccati}.

%
%
%


\section{Rates of a temporal discretization for Problem (MF-LQ)}\label{rate}
The goal of this section is to present a discretization for Problem (MF-LQ), deduce the convergence rate of this 
discretization, and propose a corresponding algorithm.

Before proving the main results, we need make some preparations. 
In Section \ref{error-of-LQ}, we have introduced two 
auxiliary stochastic LQ problems. Next we need two deterministic LQ problems.\\
\no {\bf Problem (LQ)$^{2,t}$}. For given $x_{t}\in\dbR^n$,  minimize 
\beq
\wt\cJ(t,x_t;u(\cd))
=\int_t^T\Big[\lan \wt Q(s) x(s),x(s)\ran+\lan \wt R(s) u(s),u(s)\ran \Big]\rd s+\lan \wt G x(T),x(T)\ran
\eeq
over $L^2(t,T;\dbR^m)$ subject to 
\bel{ode1}
\lt\{
\bal
&\rd x(s)
=\big[\wt A(s) x(s)+\wt B u(t)\big]\rd t\qq \forall\,s \in [t,T]\,\\
&x(t)=x_t\,,
\eal
\rt.
\ee
where
\bel{w929e7}
\lt\{
\bal
&\wt Q(\cd)=\h Q +\h C^\top P(\cd)\h C-\h C^\top P(\cd)\h D \Si_1^{-1}(\cd)\h D^\top P(\cd)\h C\,,\\
&\wt R(\cd)=\Si_1(\cd)\,,\\
&\wt G=\h G\,,\\
&\wt A(\cd)=\h A-\h B \Si_1^{-1}(\cd)\h D^\top P(\cd)\h C\,,\\
&\wt B=\h B\,,
\eal
\rt.
\ee
$P(\cd)$ solves differential Riccati equation \rf{Riccati-1}, and $\Si_1(\cd)$ is defined in \rf{w929e1}.\\
\no {\bf Problem (LQ)$_\t^{2,t_l}$}. For given $x_{t_l}\in\dbR^n$,  minimize 
\beq
\wt\cJ_\t(t_l,x_{t_l};u_\t(\cd))=\t\sum_{k=l}^{N-1}\Big[\lan\wt Q_k x_\t(t_k),x_\t(t_k)\ran+\lan \wt R_ku_\t(t_k),u_\t(t_k)\ran\Big] +\lan \wt G x_\t(T),x_\t(T)\ran
\eeq
over $\bar \dbU_\t (t_l,T)$ subject to 
\bel{odde1}
\lt\{
\bal
&x_\t(t_{k+1})=\big(I_n+\wt A_k\t\big)x_\t(t_k)+\wt B_k \t u_\t(t_k) \qq \forall\,k =l,l+1,\cds,N-1,\\
&x_\t(t_l)=x_{t_l}\,,
\eal
\rt.
\ee
where 
\bel{w929e8}
\lt\{
\bal
&\wt Q_k=\h Q +\h C^\top P_k \h C-\h C^\top P_k\h D \big(\h R+\h D^\top P_k \h D\big)^{-1}\h D^\top P_k\h C\,,\\ 
&\wt R_k=\h R+\h D^\top P_k \h D\,,\\
&\wt A_k=\h A-\h B \big(\h R+\h D^\top P_k \h D\big)^{-1}\h D^\top P_k\h C\,,\\
&\wt B_k=\h B\,,
\eal
\rt.
\ee
and 
$P_\cd$ is the solution to difference Riccati equation \rf{Riccati-D1}.

In Section \ref{error-of-LQ}, we know that Problem (LQ)$_\t^{1,t_l}$ is an approximation of Problem (LQ)$^{1,t}$. In the
below, we will prove that Problem (LQ)$^{2,t_l}_\t$ is an approximation of Problem (LQ)$^{2,t}$. These four LQ problems
are crucial in discretizing Problem (MF-LQ) and deriving the convergence rates for the corresponding discretization.
\ms

By LQ theory, to guarantee the solvability of LQ problems, one usually needs the following sufficient condition: weight coefficients are positive semi-definite/definite.
The following result is on the weight coefficients's positive semi-definiteness/definiteness of Problem (LQ)$^{2,t}$ and Problem (LQ)$_\t^{2,t_l}$.
\bl{w917l1}
{\rm(i)} For any $t\in[0,T]$, $\wt Q(t)$ is positive semi-definite and $\wt R(t)$ is uniformly positive definite, where $\wt Q(\cd)$
and $\wt R(\cd)$ are defined in \rf{w929e7}.

{\rm(ii)} For any $k=0,1,\cds,N-1$, $\wt Q_k$ is positive semi-definite and $\wt R_k$ is uniformly positive definite, where
$\wt Q_\cd$ and $\wt R_\cd$ are given in \rf{w929e8}. 
\el

\begin{proof}
We only prove assertion (i), and (ii) can be deduced in the same vein. 

By Lemma  \ref{Yong-Riccati}, we know that $P(\cd)$ is positive semi-definite, and then
\beq
\wt R(t)=\h R+\h D^\top P(t)\h D\geq \h R>0 \qq \forall\,t\in[0,T]\,,
\eeq
which is the second assertion.

To prove the first assertion, we need the following claim: For any $K\in \dbR^{n\times m}, R_0\in \dbR^{m\times m} $
satisfying $R_0>0$, it holds that
\bel{inver1}
\|K(R_0+K^\top K )^{-1}K^\top\|\leq 1\,.
\ee
Indeed, by setting $K_0=KR_0^{-1/2}$, we can deduce that
\beq
K(R_0+K^\top K )^{-1}K^\top=K_0(I+K_0^\top K_0 )^{-1}K_0^\top\,.
\eeq
Without loss of generality, we assume that $n\leq m$. By singular value decomposition, it follows that
\beq
K_0=U\begin{pmatrix}\L & 0\end{pmatrix} V^\top\,,
\eeq
where $U\in \dbR^{n\times n}$, $V\in \dbR^{m\times m}$ are orthogonal matrices, and  
$\L=\mbox{diag}\{\l_1,\cds,\l_n\}$. Hence,
\beq
\|K(R_0+K^\top K )^{-1}K^\top\|
&=\Big\| \begin{pmatrix}\L & 0\end{pmatrix} \Big(I_m+\begin{pmatrix}\L & 0\end{pmatrix}^\top \begin{pmatrix}\L & 0\end{pmatrix}\Big)^{-1}\begin{pmatrix}\L & 0\end{pmatrix}^\top \Big\|\\
&=\|\L(I_n+\L^2)^{-1}\L\|\\
&=\max_{i=1,2,\cds,n}\frac{\l_i^2}{1+\l_i^2}\\
&\leq 1\,,
\eeq
that prove \rf{inver1}.

For any $t\in[0,T]$, by taking $K=\sqrt{P(t)}D,\, R_0=\h R$ in \rf{inver1}, we have
\beq
\lt\|\big(\sqrt{P(t)}D\big)\Si_1^{-1}(t)\big(\sqrt{P(t)}D\big)^\top\rt\|
\leq 1\,.
\eeq
Subsequently, for any $x\in \dbR^n$, $t\in[0,T]$,
\beq
&\lt\lan \big[P(t)-P(t)\h D\Si_1^{-1}(t)\h D^\top P(t)\big]x,x\rt\ran\\
&=\lt\| \sqrt{P(t)}x\rt\|^2- \llan \big[ \big(\sqrt{P(t)}\h D\big)\Si_1^{-1}(t)\big(\sqrt{P(t)}\h D\big)^\top\big]\sqrt{P(t)}x, \sqrt{P(t)}x \rran\\
&\geq 0\,.
\eeq

%
Finally, by the assumption {\rm(A)}, we arrive at
\beq
\wt Q(t)=\h Q +\h C^\top \big[P(t)- P(t)\h D \Si_1^{-1}(t)\h D^\top P(t)\big]\h C\geq \h Q\geq 0 \qq \forall\, t\in[0,T]\,.
\eeq
That completes the proof.
\end{proof}


By Lemma \ref{w917l1} (i), we know that for any $t\in [0,T]$, Problem (LQ)$^{2,t}$ is uniquely solvable. Furthermore, it is easy to check
that the corresponding Riccati equation is just \rf{Riccati-2}.

In the below, we introduce a difference Riccati equation which is utilized to provide a
control's state feedback form for Problem (LQ)$_\t^{2,t_l}$:
\bel{Riccati-D3}
\lt\{
\bal
&\Pi_k=(I_n+\wt A_k\t)^\top \Pi_{k+1}(I_n+\wt A_k\t)+\t \wt Q_k-\t H_k^\top W_k^{-1}H_k \qq k=0,1,\cds,N-1,\\
&\Pi_N=\wt  G\,,\\
&H_k=\wt B_k^\top \Pi_{k+1}(I_n+\wt A_k\t)\,,\q
W_k=\wt R_k+ \t \wt B_k^\top \Pi_{k+1}\wt B_k\,,
\eal
\rt.
\ee
where $\wt A_\cd,\,\wt Q_\cd,\,\wt G,\,\wt B_\cd,\, \wt R_\cd$ are defined in \rf{w929e7} and \rf{w929e8}.

\bl{w229l2}
Let $\Pi(\cd),\,\Pi_\cd$ be the solutions to Riccati equations \rf{Riccati-2}, \rf{Riccati-D3} respectively.
Then there exists a constant $\cC$ 
independent of $\t$ such that
\bel{w917e1}
\lt\{
\bal
& \sup_{t\in[0,T]}\|\Pi(t)\| +\max_{l=0,1,\cds,N}\|\Pi_l\| \leq \cC\,, \\
& \|\Pi(t)-\Pi(s)\| \leq \cC|t-s|\qq \forall\, t,s\in [0,T]\,.
\eal
\rt.
\ee
\el

\begin{proof}
Suppose that  the optimal pairs of Problem (LQ)$^{2,t}$ and Problem (LQ)$_\t^{2,t_l}$ are $(x^*(\cd),u^*(\cd))$ and
$(x^*_\t(\cd), u^*_\t(\cd))$ respectively.

We know that 
\bel{w917e8}
\bal
\inf_{u(\cd)\in L^2(t,T;\dbR^m)}\wt \cJ(t,x_t;u(\cd))&=\lan \Pi(t)x_t,x_t\ran \qq \forall\,t\in [0,T]\,,\\
\inf_{u_\t(\cd)\in \bar \dbU_\t(t_l,T)}\wt \cJ_\t(t_l,x_{t_l};u_\t(\cd))&=\lan \Pi_lx_{t_l},x_{t_l}\ran\qq\forall\,l=0,1,\cds,N\,.
\eal
\ee

By taking $u(\cd)\equiv 0,\,u_\t(\cd)\equiv0$ in \rf{ode1}, \rf{odde1} respectively, we can deduce that
\beq
\sup_{s\in[t,T]}\|x(s)\|\leq \cC \|x_t\|\qq
\max_{k=l,l+1,\cds,N}\|x_\t(t_k)\|\leq \cC\|x_{t_l}\|\,,
\eeq
which, together with \rf{w917e8}, \rf{w301e7}$_1$ and the fact that $\|\Si_1^{-1}(\cd)\|\leq \|\h R^{-1}\|$,
$\|\big(\h R+\h D^\top P_k \h D\big)^{-1}\|\leq \|\h R^{-1}\|$, yields
\beq
&\lan \Pi(t)x_t,x_t\ran\leq \wt \cJ(t,x_t;0)=\int_t^T \lan \wt Q(s) x(s),x(s)\ran \rd s+\lan \wt G x(T),x(T)\ran
\leq \cC\|x_t\|^2\,,\\
&\lan \Pi_lx_{t_l},x_{t_l}\ran\leq \wt \cJ_\t(t_l,x_{t_l};0)
=\t\sum_{k=l}^{N-1}\lan\wt Q_k x_\t(t_k),x_\t(t_k)\ran+\lan \wt G x_\t(T),x_\t(T)\ran
\leq \cC\|x_{t_l}\|^2\,.
\eeq
Then, noting that $\Pi(\cd),\,\Pi_\cd$ are positive semi-definite, we can derive \rf{w917e1}$_1$ by the above two estimate.
Subsequently, by utilizing \rf{w301e7}$_1$, we can derive \rf{w917e1}$_2$.
\end{proof}

\bl{w922l1}
For any $t\in [0,T]$, suppose that $(x^*(\cd),\,u^*(\cd))$ is the optimal pair of Problem \rm{(LQ)$^{2,t}$}. Then, there exists
a constant $\cC$ independent of $t$ such that
\bel{w922e1}
\lt\{
\bal
&\sup_{s\in [t,T]}\big[\|x^*(s)\|+\|u^*(s)\|\big]\leq \cC \|x_t\|\,,\\
&\|x^*(s_1)-x^*(s_2)\|+\|u^*(s_1)-u^*(s_2)\|\leq \cC|s_1-s_2|\|x_t\| \qq \forall\, s_1,\,s_2\in [t,T]\,.
\eal
\rt.
\ee
\el

\begin{proof}
The desired estimate can be deduced by the optimal control's feedback representation of Problem \rm{(LQ)$^{2,t}$} 
(see, \cite[Chapter 6.2]{Yong-Zhou99})
\beq
u^*(\cd)=-\wt R^{-1}(\cd) \wt B^\top \Pi(\cd) x^*(\cd)\,,
\eeq
and Lemmas \ref{w917l1}, \ref{w229l2}.
\end{proof}

\bl{w922l3}
For any $l=0,1,\cds,N-1$, suppose that $(x_\t^*(\cd),\,u_\t^*(\cd))$ is the optimal pair of Problem \rm{(LQ)$_\t^{2,t_l}$}. Then, there exists
a positive constant $\cC$ independent of $l$ such that
\beq
\bal
\sup_{k=l,l+1,\cds,N-1}\big[\|x_\t^*(t_k)\|+\|u_\t^*(t_k)\|\big]\leq \cC \|x_{t_l}\|\,.
\eal
\eeq
\el

\begin{proof}
The desired result can be deduced by the optimal control's feedback representation  
(see, e.g., \cite[Theorem 4.3]{Zhou02})
\beq
u_\t^*(t_k)=-\big(\wt R_k+\t\wt B_k^\top \Pi_{k+1}\wt B_k\big)^{-1}\wt B_k^\top \Pi_{k+1}\big(I_n+\wt A_k\t\big) x_\t^*(t_k)
\qq k=l,l+1,\cds,N-1\,,
\eeq
and Lemmas \ref{w917l1}, \ref{w229l2}.
\end{proof}

\bl{w918l1}
Suppose that $u^*(\cd)$ is the optimal control of Problem \rm{(LQ)$^{2,t_l}$}, and $\bar x_\t(\cd)$ is the solution to
\rf{odde1} with $u_\t(\cd)=u^*(\nu(\cd))$, where $\nu(\cd)$ is defined in \rf{w128e1}. Then, it holds that
\bel{w918e2}
\lt\{
\bal
&\max_{k=l,l+1,\cds,N}\|\bar x_\t(t_k)\|\leq \cC \|x_{t_l}\|\,,\\
& \max_{k=l,l+1,\cds,N}\|\bar x_\t(t_k)-x^*(t_k)\|\leq \cC\t\|x_{t_l}\|\,.
\eal
\rt.
\ee

\el

\begin{proof}
By \rf{odde1}, we can derive
\beq
\bal
\bar x_\t(t_{k+1})&=\big(I_n+\wt A_k\t\big)\bar x_\t(t_k)+\wt B_k \t u^*(t_k)\\
&=\big(I_n+\wt A_k\t\big)\big[\big(I_n+\wt A_{k-1}\t\big)\bar x_\t(t_{k-1})+\wt B_{k-1} \t u^*(t_{k-1})\big]+\wt B_k \t u^*(t_k)\\
&=\cds\\
&=\prod_{i=k}^l \big(I_n+\wt A_i\t\big)x_{t_l}+\sum_{j=l}^k\prod_{i=k}^{j+1}\big(I_n+\wt A_i\t\big)\wt B_j \t u^*(t_{j})\,,
\eal
\eeq
which, together with \rf{w922e1}$_1$ and the fact that
\beq
\max_{i=l,l+1,\cds,N-1}\|I_n+\wt A_i \t\|\leq 1+\cC \t\,,\qq \max_{j=l,l+1,\cds,N-1}\|\wt B_j\|=\|\h B\|\leq \cC\,,
\eeq
leads to \rf{w918e2}$_1$.

By setting $e_k=\bar x_\t(t_k)-x^*(t_k)$, we can arrive at
\beq
&e_{k+1}-e_k\\
&=\int_{t_k}^{t_{k+1}} \Big[\big[\wt A_k-\wt A(t)\big] \bar x_\t(t_k)
+\wt A(t) e_k
+\wt A(t) \big[ x^*(t_k)-x^*(t) \big]
+\wt B\big[ u^*(t_k)-u^*(t)\big] \Big]\rd t\,.
\eeq
Subsequently, by applying Lemma \ref{w922l1} and equation \rf{odde1}, we find that
\beq
\|e_{k+1}\|\leq (1+\cC\t)\|e_k\|+\cC \t^2\,,
\eeq
which, together with discrete Gronwall's inequality, yields \rf{w918e2}$_2$.
\end{proof}

The following result can be deduced by the similar trick used in the proof of Lemma \ref{w918l1}.
\bl{w922l4}
Suppose that $u_\t^*(\cd)$ is the optimal control of Problem \rm{(LQ)$_\t^{2,t_l}$}, and $\check x(\cd)$ is the solution to
\rf{ode1} with $u(\cd)=u_\t^*(\cd)$. Then, it holds that
\beq
\lt\{
\bal
&\sup_{t\in [t_l,T]}\|\check x(t)\|\leq \cC \|x_{t_l}\|\,,\\
& \max_{k=l,l+1,\cds,N}\|\check x(t_k)-x_\t^*(t_k)\|\leq \cC\t\|x_{t_l}\|\,.
\eal
\rt.
\eeq

\el

\bl{w922l2}
Let $\Pi(\cd),\,\Pi_\cd$ be the solutions to Riccati equations \eqref{Riccati-2}, \rf{Riccati-D3} respectively.
Then there exists a positive constant $\cC$ 
independent of $\t$ such that
\bel{w922e3}
\max_{l=0,1,\cds,N}\|\Pi(t_l)-\Pi_l\|\leq \cC \t\,.
\ee
\el

\begin{proof}
Suppose that  the optimal pairs of Problem (LQ)$^{2,t_l}$ and Problem (LQ)$_\t^{2,t_l}$ are $(x^*(\cd),u^*(\cd))$ and
$(x^*_\t(\cd), u^*_\t(\cd))$ respectively. 
For any $l=0,1,\cds,N-1$, $x_{t_l}\in \dbR^n$, we consider the following two cases. 

\no{\bf Case I:} $\lan \Pi(t_l)x_{t_l}, x_{t_l}\ran\leq \lan \Pi_l x_{t_l},x_{t_l}\ran$.

In this case, we have
\beq
\bal
&\lan \Pi_l x_{t_l}, x_{t_l}\ran - \lan \Pi(t_l) x_{t_l},x_{t_l}\ran =\wt\cJ_\t(t_l,x_{t_l};u_\t^*(\cd))-\wt\cJ(t_l,x_{t_l};u^*(\cd))\\
&\leq \wt\cJ_\t(t_l,x_{t_l};u^*(\nu(\cd)))-\wt\cJ(t_l,x_{t_l};u^*(\cd))\\
&=\sum_{k=l}^{N-1}\int_{t_k}^{t_{k+1}}\Big[\lan\wt Q_k \bar x_\t(t_k),\bar x_\t(t_k)\ran-\lan\wt Q(t) x^*(t),x^*(t)\ran\\
&\qq\qq\qq+\lan \wt R_ku^*(t_k),u^*(t_k)\ran-\lan\wt R u^*(t),u^*(t)\ran\Big] \rd t\\
&\q+\lan \wt G \bar x_\t(T),\bar x_\t(T)\ran-\lan \wt G x^*(T),x^*(T)\ran\\
&\leq \sum_{k=l}^{N-1}\int_{t_k}^{t_{k+1}}
\Big[  \|\wt Q_k-\wt Q(t)\|\|\bar x_\t(t_k)\|^2 +\big[\|\wt Q(t) \bar x_\t(t_k)\|+\|\wt Q(t) x^*(t)\|\big]\|\bar x_\t(t_k)-x^*(t)\|\\
&\q+\|\wt R_k-\wt R(t)\|\|\bar u^*(t_k)\|^2 +\big[\|\wt R(t) u^*(t_k)\|+\|\wt R(t) u^*(t)\|\big]\|u^*(t_k)-u^*(t)\|\Big] \rd t\\
&\q+\big[\|\wt G \bar x_\t(T)\|+\|\wt G x^*(T)\|\big]\|\bar x_\t(T)-x^*(T)\|\\
&\leq \cC\t\|x_{t_l}\|^2\,.
\eal
\eeq
Here, we utilize Lemmas \ref{w229l4}, \ref{w229l2}, \ref{w922l1} and \ref{w918l1}.

\no{\bf Case II:} $\lan \Pi(t_l)x_{t_l}, x_{t_l}\ran> \lan \Pi_l x_{t_l},x_{t_l}\ran$.

In this case, we can deduce that
\beq
\bal
&\lan \Pi(t_l) x_{t_l},x_{t_l}\ran-\lan \Pi_l x_{t_l}, x_{t_l}\ran  \leq\wt\cJ(t_l,x_{t_l};u_\t^*(\cd))-\wt\cJ(t_l,x_{t_l};u_\t^*(\cd))\\
&= \sum_{k=l}^{N-1}\int_{t_k}^{t_{k+1}}\Big[\lan\wt Q(t) \check x(t),\check x(t)\ran-\lan\wt Q_k x_\t^*(t_k),x_\t^*(t_k)\ran+\lan \big(\wt R(t)-\wt R_k\big) u_\t^*(t_k),u_\t^*(t_k)\ran\Big] \rd t\\
&\q+\lan \wt G \check x(T),\check x(T)\ran-\lan \wt G x_\t^*(T),x_\t^*(T)\ran\\
&\leq \sum_{k=l}^{N-1}\int_{t_k}^{t_{k+1}}
\Big[  \|\wt Q_k-\wt Q(t)\|\| x_\t^*(t_k)\|^2 +\big[\|\wt Q(t) \check x(t)\|+\|\wt Q(t) x_\t^*(t_k)\|\big]\|x_\t^*(t_k)-\check x(t)\|\\
&\q+\|\wt R_k-\wt R(t)\|\| u_\t^*(t_k)\|^2 \Big] \rd t
+\big[\|\wt G \check x(T)\|+\|\wt G x_\t^*(T)\|\big]\|\check x(T)-x_\t^*(T)\|\\
&\leq \cC\t\|x_{t_l}\|^2\,.
\eal
\eeq
Here, we apply Lemmas \ref{w229l4}, \ref{w229l2}, \ref{w922l3} and \ref{w922l4}.

Finally, for any symmetric $K$, it holds
\bel{norm-K}
\|K\|=\max_{\|x\|=1}|\lan Kx,x \ran|\,.
\ee
Combining with the above two cases and \rf{norm-K}, we can derive \rf{w922e3}. That completes the proof.
\end{proof}

The following is the main theorem.

\bt{rate-tx}
Under the assumption {\rm (A)},
suppose that $(x^*(\cd),u^*(\cd))$  is the optimal pair of Problem {\rm (MF-LQ)}, and $x_\t(\cd)$ is the state to system
\rf{mf-sdde} with a closed--loop control:
\bel{w923e2}
u_\t(t_k)=-(W_k^1)^{-1}H_k^1(x_\t(t_k)-\me[x_\t(t_k)])-(W_k^2)^{-1}H_k^2\me[x_\t(t_k)] \qq k=0,1,\cds,N-1\,,
\ee
where $W_\cd^1,\,H_\cd^1$ are defined in \rf{w929e4}, 
\bel{w923e3}
\bal
W_k^2=\h R+\h B^\top\Pi_{k+1}\h B\t+\h DP_{k+1}\h D\,,\qq
H_k^2=\h B^\top \Pi_{k+1}(I_n+\h A\t)+\h D^\top P_{k+1}\h C\,,
\eal
\ee
and $P_\cd,\,\Pi_\cd$ solve difference Riccati equations \rf{Riccati-D1}, \rf{Riccati-D3} respectively.
Then there exists a constant $\cC$ independent of $\t$ such that
\bel{rates}
\lt\{
\bal
&\sup_{k=0,1,\cds,N}\big[\|\me[x^*(t_k)]-\me[x_\t(t_k)]\|+\|\me[u^*(t_k)]-\me[u_\t(t_k)]\| \big] \leq \cC \sqrt{\t}\|x_0\|\,,\\
&\sup_{k=0,1,\cds,N}\me\big[\|x^*(t_k)-x_\t(t_k)\|^2+\|u^*(t_k)-u_\t(t_k)\|^2 \big] \leq \cC \t\|x_0\|^2\,.
\eal
\rt.
\ee
Furthermore, if $C=\bar C=0$, then the following sharp estimate holds
\bel{rates2}
\bal
\sup_{k=0,1,\cds,N}\big[\|\me[x^*(t_k)]-\me[x_\t(t_k)]\|+\|\me[u^*(t_k)]-\me[u_\t(t_k)]\| \big] \leq \cC \t\|x_0\|\,.
\eal
\ee
\et

\begin{proof}
The proof is long, and we divide it into two steps.

\no{\bf (1).}  By Lemma \ref{Yong-Riccati}, we take the optimal control's state feedback form \rf{w916e1} for Problem (MF-LQ), and then $\me[x^*(\cd)]$
satisfies
\beq
\lt\{
\bal
&\rd \me[x^*(t)]=\big[\h A-\h B \Si_1^{-1}(t)\big(\h B^\top \Pi(t) +\h D^\top P(t)\h C\big)\big]\me[x^*(t)]\rd t \qq \forall\,t\in [0,T]\,,\\
&\me[x^*(0)]=x_0\,.
\eal
\rt.
\eeq 
For system \rf{mf-sdde}, we use the control's state feedback form \rf{w923e2} which may be not optimal, and
derive the following difference equation
\beq
\lt\{
\bal
& \me[x_\t(t_{k+1})]=\big[I_n+\h A\t-\t\h B (W_k^2)^{-1}H_k^2 \big]\me[x_\t(t_k)] \qq k=0,1,\cds,N-1\,,\\
&\me[x_\t(0)]=x_0\,.
\eal
\rt.
\eeq 
It is easy to find that
\bel{w917e3}
\begin{array}{c}
\ds\max_{k=0,1,\cds,N-1}\Big[\sup_{t\in[t_k,t_{k+1})}\lt\|\me[x^*(t)]-\me[x^*(t_k)]\rt\|
+\|\me[x_\t(t_{k+1})]-\me[x_\t(t_k)]\|\Big]\leq \cC\t\,,\\
\ds \sup_{t\in [0,T]}\|\me[x^*(t)]\|+\max_{k=0,1,\cds,N}\|\me[x_\t(t_k)]\|\leq \cC\,.
\end{array}
\ee

By setting $\check e_k=\me[x^*(t_k)]-\me[x_\t(t_k)]$, we can arrive at
\bel{w917e2}
\bal
\check e_{k+1}
&=\check e_k+\int_{t_k}^{t_{k+1}}\Big\{ \big[\h A-\h B \Si_1^{-1}(t)\big(\h B^\top \Pi(t) +\h D^\top P(t)\h C\big)\big]\big[\me[x^*(t)]-\me[x^*(t_k)]\big]\\
&\q +\big[\h A-\h B \Si_1^{-1}(t)\big(\h B^\top \Pi(t) +\h D^\top P(t)\h C\big)\big] \check e_k\\
&\q-\h B \Big[  \Si_1^{-1}(t)\big(\h B^\top \Pi(t) +\h D^\top P(t)\h C\big) -(W_k^2)^{-1}H_k^2\Big]  \me[x_\t(t_k)]     \Big\}\rd t\,.
\eal
\ee
Now, applying Lemmas \ref{w229l4}, \ref{w229l2} and the assumption \rm{(A)}, we can deduce that
\bel{w922e11}
\bal
\sup_{t\in[0,T]}\lt\|\h A-\h B \Si_1^{-1}(t)\big(\h B^\top \Pi(t) +\h D^\top P(t)\h C\big)\rt\|\leq \cC\,.
\eal
\ee
By \rf{w301e7}$_2$, \rf{w917e1}$_1$, and the fact that $P(\cd),\,P_\cd,\,\Pi_\cd$ are positive semi-definite, we have the following estimate: for any $t\in[t_k,t_{k+1}]$,
\beq
\lt\|\Si_1^{-1}(t)-(W_k^2)^{-1}\rt\|
&=\lt\|\Si_1^{-1}(t)\big[W_k^2-\Si_1(t)\big]
(W_k^2)^{-1}\rt\|\\
&=\lt\|\Si_1^{-1}(t)\Big[\h D^\top \big(P_{k+1}-P(t)\big)\h D+\h B^\top \Pi_{k+1}\h B\t\Big]
(W_k^2)^{-1}\rt\|\\
&\leq \cC\t \|\h R^{-1}\|^2\,,
\eeq
which, together with \rf{w301e7}$_3$ in Lemma \ref{w229l4} and Lemma \ref{w922l2},  leads to
\bel{w922e12}
\lt\|  \Si_1^{-1}(t)\big(\h B^\top \Pi(t) +\h D^\top P(t)\h C\big) -(W_k^2)^{-1}H_k^2\rt\|\leq \cC\sqrt{\t} \qq \forall\,t\in[0,T]\,.
\ee
Subsequently, by virtue of  \rf{w917e3}, \rf{w917e2}, \rf{w922e11}, \rf{w922e12} and discrete Gronwall's inequality, we 
can derive 
\beq
\max_{k=0,1,\cds,N}\|\check e_k\|\leq \cC\sqrt{\t}\,.
\eeq

Next, applying optimal controls' feedback representations \rf{w916e1} and \rf{w923e2}, we can conclude that
\beq
\|\me[u^*(t_k)]-\me[u_\t(t_k)]\|\leq \cC\sqrt{\t}\,.
\eeq
Hence we get \rf{rates}$_1$.

When $C=\bar C=0$, \rf{w922e12} turns to
\beq
\lt\|  \Si_1^{-1}(t)\big(\h B^\top \Pi(t) +\h D^\top P(t)\h C\big) -(W_k^2)^{-1}H_k^2\rt\|\leq \cC\t \qq \forall\,t\in[0,T]\,,
\eeq
which can be used to derive \rf{rates2}.

\ms

\no{\bf (2).}  In this step we tend to prove the estimate \rf{rates}$_2$.
Based on the feedback  representations \rf{w916e1} and \rf{w923e2}, the state equations 
\rf{mf-sde} (with $t$ suppresed) and 
\rf{mf-sdde} satisfy
\beq
\lt\{
\bal
&\rd x^*=\Big\{A x^*-B\Big[\Si_0^{-1}\big(B^\top P+D^\top PC\big)\big(x^*-\me[x^*]\big)
+\Si_1^{-1}\big(\h B^\top \Pi +\h D^\top P\h C\big)\me[x^*] \Big]\\
&\qq\qq+\bar A \me[x^*]+\bar B\me[u^*] \Big\}\rd t\\
&\q\qq+\Big\{Cx^*-D\Big[\Si_0^{-1}\big(B^\top P+D^\top PC\big)\big(x^*-\me[x^*]\big)
+\Si_1^{-1}\big(\h B^\top \Pi +\h D^\top P\h C\big)\me[x^*] \Big]\\
&\qq\qq+\bar C \me[x^*]+\bar D\me[u^*] \Big\}\rd W(t) \qq \forall\,t\in [0,T]\,,\\
&x^*(0)=x_0\,,
\eal
\rt.
\eeq 
\beq
\lt\{
\bal
& x_\t(t_{k+1})=\big[I_n+A\t\big]x_\t(t_k)+\t\bar A\me[x_\t(t_k)]+\t\bar B\me[u_\t(t_k)]\\
&\qq\qq-\t B\big[(W_k^1)^{-1}H_k^1(x_\t(t_k)-\me[x_\t(t_k)])+(W_k^2)^{-1}H_k^2\me[x_\t(t_k)]\big] \\
&\qq\qq+\Big[ Cx_\t(t_k)+ \bar C\me[x_\t(t_k)] +\bar D\me[u_\t(t_k)]\\
&\qq\qq- D\big[(W_k^1)^{-1}H_k^1(x_\t(t_k)-\me[x_\t(t_k)])+(W_k^2)^{-1}H_k^2\me[x_\t(t_k)]\big]\Big]\D_{k+1}W  \\
&\qq\qq\qq\qq\qq\qq\forall\, k=0,1,\cds,N-1    \,,\\
&x_\t(0)=x_0\,.
\eal
\rt.
\eeq 

Set $\bar e_k=x^*(t_k)-x_\t(t_k)$. Then by above two equations, we can see that
\begin{eqnarray}
&&\bar e_{k+1}-\bar e_k \nonumber\\
&&=\int_{t_k}^{t_{k+1}}\Big\{A(x^*(t)-x_\t(t_k))+\bar A\big(\me[x^*(t)]-\me[x_\t(t_k)]\big) +\bar B \big(\me[u^*(t)]-\me[u_\t(t_k)]\big)  \nonumber\\
&&\q-B \Big[ \Si_0^{-1}\big(B^\top P(t)+D^\top P(t)C\big)\big(x^*(t)-\me[x^*(t)]\big)  -(W_k^1)^{-1}H_k^1(x_\t(t_k)-\me[x_\t(t_k)])\Big] \nonumber\\
&&\q-B\Big[ \Si_1^{-1}\big(\h B^\top \Pi(t) +\h D^\top P(t)\h C\big)\me[x^*(t)]-(W_k^2)^{-1}H_k^2\me[x_\t(t_k)]\Big]  \Big\}\rd t \nonumber\\
&&\q+\int_{t_k}^{t_{k+1}}\Big\{C(x^*(t)-x_\t(t_k))+\bar C\big(\me[x^*(t)]-\me[x_\t(t_k)]\big)+\bar D \big(\me[u^*(t)]-\me[u_\t(t_k)]\big) \nonumber  \\
&&\q-D \Big[ \Si_0^{-1}\big(B^\top P(t)+D^\top P(t)C\big)\big(x^*(t)-\me[x^*(t)]\big)  -(W_k^1)^{-1}H_k^1(x_\t(t_k)-\me[x_\t(t_k)])\Big] \nonumber \\
&&\q-D\Big[ \Si_1^{-1}\big(\h B^\top \Pi(t) +\h D^\top P(t)\h C\big)\me[x^*(t)]-(W_k^2)^{-1}H_k^2\me[x_\t(t_k)]\Big]  \Big\}\rd W(t) \nonumber\\
&&:=\int_{t_k}^{t_{k+1}} I(t)\rd t+\int_{t_k}^{t_{k+1}} J(t)\rd W(t)\,. \label{w917e4}
\end{eqnarray}
By multiplying $\bar e_{k+1}$, and then taking conditional expectation on both sides of \rf{w917e4}, we can find that
\bel{w917e5}
\mbox{left side}=\frac 1 2 \me\big[\|\bar e_{k+1}\|^2-\|\bar e_k\|^2+\|\bar e_{k+1}-\bar e_k\|^2\big]\,.
\ee

Now we estimate the derived terms on the right side of \rf{w917e4}.
\beq
&\me\Big[\int_{t_k}^{t_{k+1}}\lan A(x^*(t)-x_\t(t_k)),\bar e_{k+1}\ran \rd t\Big]\\
&=\me\Big[\int_{t_k}^{t_{k+1}}\lan A(x^*(t)-x^*(t_k))+A\bar e_k,\bar e_{k+1}\ran \rd t\Big]\\
&\leq \frac {\|A\|}{2}\me\Big[\int_{t_k}^{t_{k+1}} \|x^*(t)-x^*(t_k)\|^2+\|\bar e_{k+1}\|^2\Big]\rd t+\frac{\|A\|}{2}\t\me\big[\|\bar e_k\|^2+\|\bar e_{k+1}\|^2\big]\\
&\leq \cC\t^2+\cC\t \me\big[\|\bar e_k\|^2+\|\bar e_{k+1}\|^2\big]\,.
\eeq
Here we apply the fact that $\me[ \|x^*(t)-x^*(s)\|^2]\leq \cC|t-s|$. In the same vein, we can estimate other Lebesgue integrals and deduce that
\bel{w917e6}
\me\Big[\int_{t_k}^{t_{k+1}} \lan I(t),\bar e_{k+1}\ran \rd t\Big]\leq \cC\t^2+\cC\t \me\big[\|\bar e_k\|^2+\|\bar e_{k+1}\|^2\big]\,.
\ee

For the It\^o integrals, following the procedure in Step (1), we can obtain
\beq
&\me\Big\lan \int_{t_k}^{t_{k+1}} D\Big[ \Si_1^{-1}\big(\h B^\top \Pi(t) +\h D^\top P(t)\h C\big) \me[x^*(t)]-(W_k^2)^{-1}H_k^2\me[x_\t^*(t_k)]\Big] \rd W(t) , \bar e_{k+1}\Big\ran\\
&=\me\Big\lan \int_{t_k}^{t_{k+1}} D\Big[ \Si_1^{-1}\big(\h B^\top \Pi(t) +\h D^\top P(t)\h C\big)\me[x^*(t)]-(W_k^2)^{-1}H_k^2\me[x_\t^*(t_k)]\Big] \rd W(t) , \bar e_{k+1}-\bar e_k\Big\ran\\
&\leq \frac 1 8 \me\big[\|\bar e_{k+1}-\bar e_k\|^2\big]\\
&\qq+2\me\Big[ \int_{t_k}^{t_{k+1}} \Big\|D\Big[ \Si_1^{-1}\big(\h B^\top \Pi(t) +\h D^\top P(t)\h C\big) \me[x^*(t)]-(W_k^2)^{-1}H_k^2\me[x_\t^*(t_k)]\Big]\Big\|^2 \rd t\Big]\\
&\leq \frac 1 8 \me\big[\|\bar e_{k+1}-\bar e_k\|^2\big]+\cC\t^2+\cC\t\me[\|\bar e_k\|^2]\,.
\eeq
Also, by the same trick, we can conclude that
\bel{w917e7}
\me\Big[\int_{t_k}^{t_{k+1}} \lan J(t),\bar e_{k+1}\ran \rd W(t)\Big]\leq \frac 3 8 \me\big[\|\bar e_{k+1}-\bar e_k\|^2\big]+\cC\t^2+\cC\t\me[\|\bar e_k\|^2]\,.
\ee

Finally, by \rf{w917e5}, \rf{w917e6} and  \rf{w917e7}, for sufficiently small $\t$, we can arrive that
\beq
\me[\|\bar e_{k+1}\|^2]\leq (1+\cC\t)\me[\|\bar e_{k}\|^2]+\cC\t^2\,,
\eeq 
and then by discrete Gronwall's inequality derive
\beq
\max_{k=0,1,\cds,N}\me[\|\bar e_{k}\|^2]\leq \cC\t\,.
\eeq
Subsequently, by controls' feedback  representations \rf{w916e1} and \rf{w923e2}, we can derive
\beq
\max_{k=0,1,\cds,N}\me[\|u^*(t_k)-u_\t(t_k)\|^2]\leq \cC\t\,.
\eeq
That completes the proof.
\end{proof}

Thanks to Theorem \ref{rate-tx}, we can present the following Riccati based algorithm.
\begin{algorithm}
\caption{ Solving Problem (MF-LQ) by the explicit Euler discretization} \label{alg2}

Fix time mesh $I_\t$ with step size $\t$.
\begin{enumerate}[1.]
\item Compute $\{P_k\}_{k=0}^N,\,\{\Pi_k\}_{k=0}^N$ via discrete Riccati equations \rf{Riccati-D1}, \rf{Riccati-D3}.

\item Compute approximating mean optimal state $\{\me[x_\t(t_k)]\}_{k=0}^N$ via \rf{mf-sdde} by taking 
\bel{w923e6}
u_\t(t_k)=-(W_k^1)^{-1}H_k^1(x_\t(t_k)-\me[x_\t(t_k)])-(W_k^2)^{-1}H_k^2\me[x_\t(t_k)] \qq k=0,1,\cds,N-1\,,
\ee
where $W_k^1,\,W_k^2,\,H_k^1,\,H_k^2$ are defined by \rf{w929e4}, \rf{w923e3}.
Then obtain approximating mean optimal control by
\beq
\me[u_\t(t_k)]=-(W_k^2)^{-1}H_k^2\me[x_\t(t_k)] \qq k=0,1,\cds,N-1\,.
\eeq

\item Compute approximating optimal state $\{x_\t(t_k)\}_{k=0}^N$ via \rf{mf-sdde} by taking a closed--loop control \rf{w923e6}, where $\{\me[x_\t(t_k)]\}_{k=0}^N,\,\{\me[u_\t(t_k)]\}_{k=0}^N$ are computed in Step 2.
Then obtain approximating  optimal control $\{u_\t(t_k)\}_{k=0}^N$ by \rf{w923e6}.

\end{enumerate}
\end{algorithm}

\br{w923r1}
\begin{enumerate}[{\rm(1)}]
\item 
By virtue of cost functional $\cJ_\t(t_l,x_{t_l};\cd)$ defined in \rf{cost-l} and the controlled 
mean-field stochastic difference equation \rf{mf-sdde}, we can introduce the following LQ problem:
\no {\bf Problem (MF-LQ)$_\t$}. For given $l=0, x_{t_l}=x_0\in\dbR^n$,  minimize $\cJ_\t(0,x_0;u_\t(\cd))$ over $\dbU_\t (0,T)$ subject to \rf{mf-sdde}.

The state feedback form \rf{w923e2} may be not the optimal control of Problem {\rm(MF-LQ)$_\t$}. But when $C=\bar C=0$,
\rf{w923e2} is optimal. The reader can refer to \cite[Theorem 3.2]{Ni13}.

\item When $\bar A=\bar C=\bar Q=\bar G=0$, $\bar B=\bar D=0$, $\bar R=0$, Problem {\rm(MF-LQ)} turns to a 
stochastic LQ problem, and the convergence rate order of the proposed discretization is $\frac 1 2$, which is 
the same as the corresponding rate in \cite{Wang21Riccati}. 
In \cite{Wang21}, discretization of same LQ problems is proposed by approximating 
forward-backward  stochastic differential equations by a gradient descent method, and the convergence rate is derived. 
However, since the gradient-descent based algorithm stems from MP, 
conditional expectations in that algorithm have to be calculated/simulated, which means that algorithm in \cite{Wang21} is
 computational expensive.

\item When $A=C=\bar C= Q= G=0$, $ B=D=\bar D=0$, $R=0$, Problem {\rm(MF-LQ)} turns to a deterministic LQ problem, and the convergence rate order of the proposed discretization is $1$.

\end{enumerate}
\er

\br{307r1}
Compared with Lemma \ref{MP-1}, when approximating optimal pair of Problem {\rm(MF-LQ)} by Algorithm \ref{alg2}, 
we need not numerically solve a coupled mean-field forward-backward  stochastic differential equation, which is not a 
trivial work, because one has to simulate conditional expectations and decouple coupled equation (see, e.g., \cite{Dunst-Prohl16,Prohl-Wang21a,Prohl-Wang22,Wang21}). Therefore, Algorithm \ref{alg2} can greatly reduce computational cost.

\er

\br{930r1}
By \cite[Section 4]{Yong13SICON}, it follows that
\beq
\lt\{
\bal
&y(\cd)=P(\cd)\big( x^*(\cd)-\me[x^*(\cd)]\big)+\Pi(\cd)\me[x^*(\cd)]\,,\\
&z(\cd)=P(\cd)\big( C\lt(x^*(\cd)-\me[x^*(\cd)]\rt)+\h C\me[x^*(\cd)]
+D\lt(u^*(\cd)-\me[u^*(\cd)]\rt)+\h D\me[u^*(\cd)]\big)\,,
\eal
\rt.
\eeq
where $(y(\cd),z(\cd))$ solves mean-field backward stochastic differential equation 
(MF-BSDE, for short) 
\rf{w1109e1}$_2$.
Also, based on Algorithm \ref{alg2}, we can obtain
$(P(\cd),\Pi(\cd),x^*(\cd),\me[x^*(\cd)],u^*(\cd),\me[u^*(\cd)])$'s approximation $(P_\cd,\Pi_\cd, x_\t(\cd),\me[x_\t(\cd)],u_\t(\cd), \me[u_\t(\cd)])$.
Hence, we can numerically solve MF-BSDE \rf{w1109e1}$_2$ as follows: \\
\no
{\rm(i)} respectively approximating $\me[y(\cd)]$ and $\me[z(\cd)]$ 
\beq
\me[y_\t(\cd)]=\Pi_{\pi(\cd)} \me[x_\t(\pi(\cd))] \q\mbox{and}\q \me[z_\t(\cd)]=P_{\pi(\cd)} \big( \h C \me[x_\t(\pi(\cd))]+\h D\me[u_\t(\pi(\cd))]\big);
\eeq
\no
{\rm(ii)} respectively approximating $y(\cd)$ and $z(\cd)$ by 
\beq
y_\t(\cd)=\Pi_{\pi(\cd)} \big(x_\t(\pi(\cd))-\me[x_\t(\pi(\cd))]\big)+\me[y_\t(\cd)]
\eeq
and
\beq
z_\t(\cd)=P_{\pi(\cd)} \big( C(x_\t(\pi(\cd))- \me[x_\t(\pi(\cd))])+ D(u_\t(\pi(\cd))-\me[u_\t(\pi(\cd))])\big)+\me[z_\t(\cd)]\,,
\eeq
where $\pi: [0,T]\to \dbN $ 
\beq
\pi (t)=k\qq \forall\, t\in [t_k,t_{k+1}),\,k=0,1,\cds,N-1\,.
\eeq
By above two steps, we actually propose an effective algorithm to solve MF-BSDE \rf{w1109e1}$_2$. Furthermore, applying 
Lemmas \ref{w301e7}, \ref{w922l2} and Theorem \ref{rate-tx}, we can derive the convergence rates for this algorithm:
\beq
\lt\{
\bal
&\sup_{k=0,1,\cds,N}\big[\|\me[y(t_k)]-\me[y_\t(t_k)]\|+\|\me[z(t_k)]-\me[z_\t(t_k)]\| \big] \leq \cC \sqrt{\t}\,,\\
&\sup_{k=0,1,\cds,N}\me\big[\|y(t_k)-y_\t(t_k)\|^2+\|z(t_k)-z_\t(t_k)\|^2 \big] \leq \cC \t\,.
\eal
\rt.
\eeq
Moreover, when $C=\bar C=0$, it holds that
\beq
\sup_{k=0,1,\cds,N}\big[\|\me[y(t_k)]-\me[y_\t(t_k)]\| \big] \leq \cC \t \,.
\eeq

\er

In the below, we provide an example to verify our theoretical results.

\begin{example}\label{wex1}
In the current example, we take $m=n=1$, $T=1$, $A=1,\,\bar A=0,\,B=0,\,\bar B=1,\,C=0,\,\bar C=0,\,D=0,\,\bar D=1,\,Q=1,\,\bar Q=-1,\,R=1,\,\bar R=-\frac 1 2,\,G=0,\,\bar G=1$, and $x_0$=1.
Under above setting, we can derive solutions to Riccati equations \rf{Riccati-1} and \rf{Riccati-2}, reading
\beq
P(t)=\frac{e^{2-2t}-1}{2} \,,
\qq \Pi(t)= \frac{1}{3-2t}e^{2-2t} \qq \forall\,t\in[0,T]\,.
\eeq
Then, by the optimal control's state feedback form \rf{w916e1}, we can derive 
\beq
\me[x^*(t)]=\frac{3-2t}{3}e^t\,,
\qq  \me[u^*(t)]=-\frac{2}{3}e^t \qq \forall\,t\in [0,T]\,,
\eeq
and then
\beq
x^*(t)=\frac{3-2t-2W(t)}{3}e^t\,,
\qq  u^*(t)=-\frac{2}{3}e^t \qq \forall\,t\in [0,T]\,.
\eeq

Numerical results by proposed Algorithm \ref{alg2}  are presented in  Figure \ref{Fig.lable}.

\vspace{-0ex}
\begin{figure}[H] 
	\centering
\includegraphics[width = 0.9\textwidth, scale=0.45]{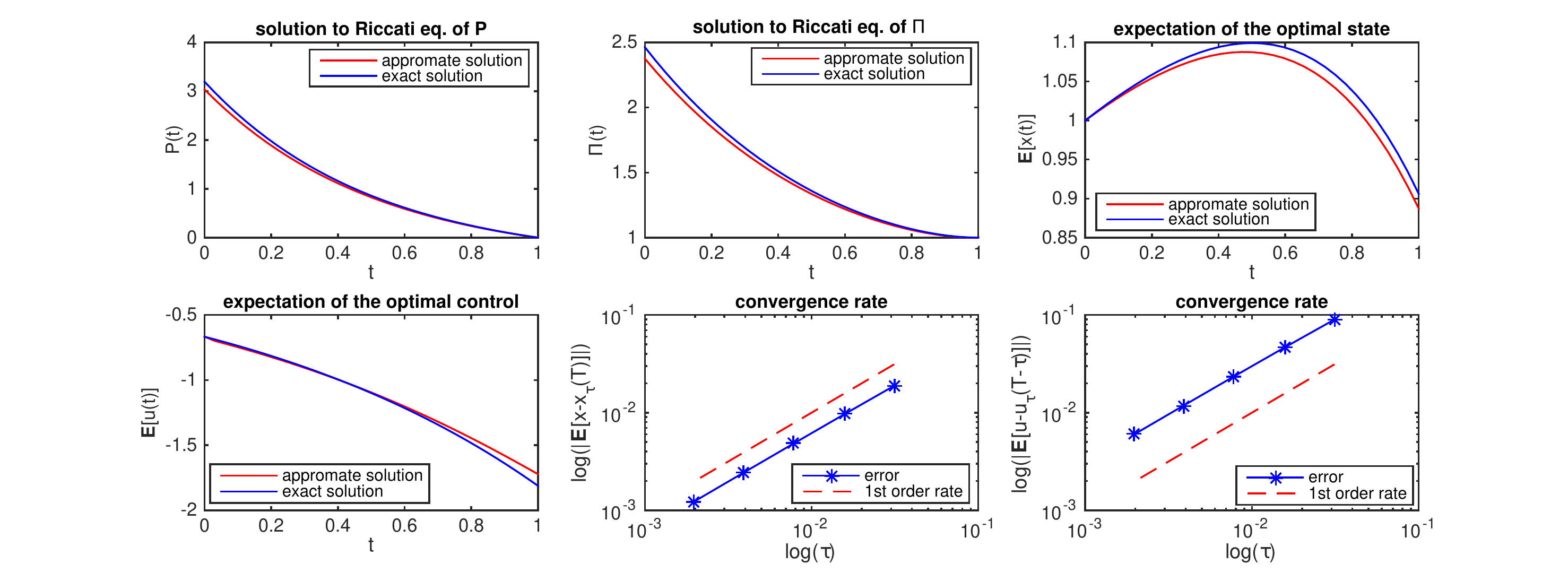}
\includegraphics[width = 0.9\textwidth, scale=0.45]{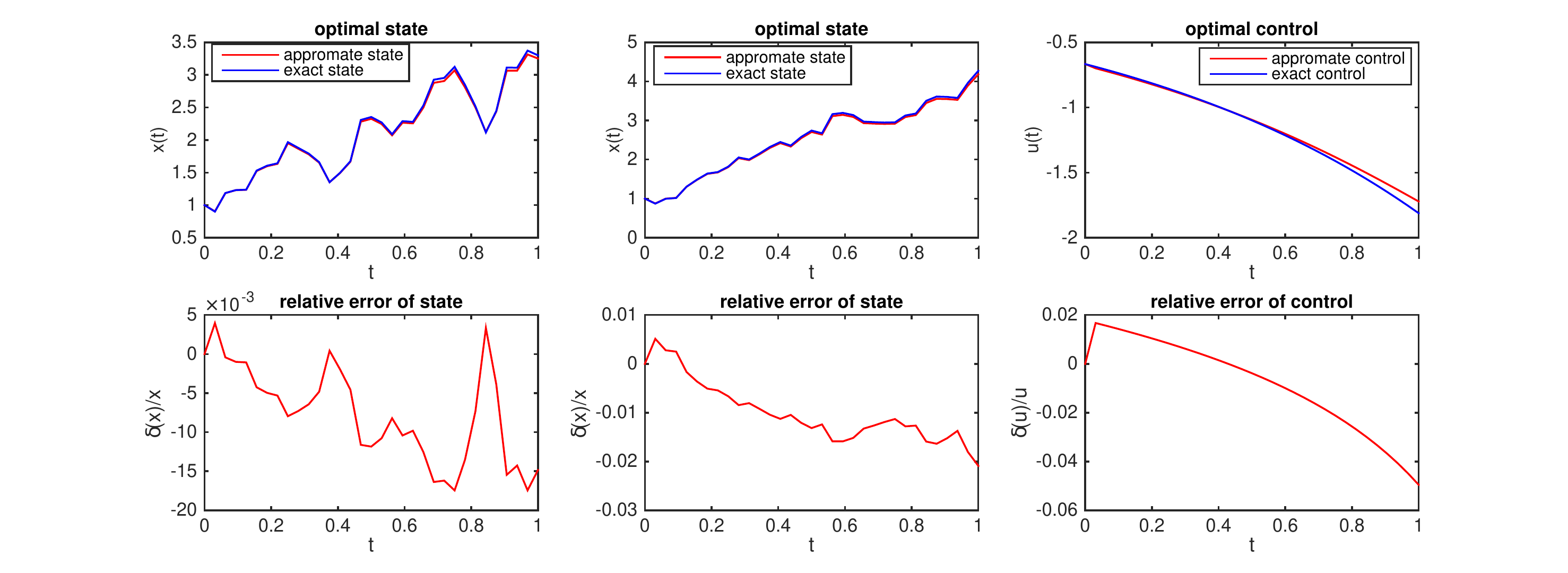}
\vspace{-3ex}
\caption{Numerical results for Example \ref{wex1}. 
$\d(x)/x:=\frac{x_\t(\cd)-x^*(\cd)}{x^*(\cd)}$ and $\d(u)/u:=\frac{u_\t(\cd)-u^*(\cd)}{u^*(\cd)}$.}
\label{Fig.lable}
\end{figure}


In the first two figures of Figure \ref{Fig.lable},  we plot the exact solutions to Riccati equations \rf{Riccati-1} and \rf{Riccati-2}, 
the numerical solutions by Algorithm \ref{alg2} respectively, with the time step $\ds \t =2^{-5}$.
In the third and fourth figures, we plot the expectation of the optimal pair, and its numerical counterpart by Algorithm \ref{alg2}, with the time step $\ds \t =2^{-5}$.  
In the fifth and sixth figures, to demonstrate the convergence rate, we  plot the corresponding convergence curves. 
For reference, we add dashed red lines of slope $1$.
In the seventh and eighth figures, we plot two path of the optimal state, and in tenth and eleventh figures, we show
the relative error of the optimal state, with the time step $\ds \t =2^{-5}$. The remaining two figures show the error and the relative error of the optimal
control, with the time step $\ds \t =2^{-5}$.

\end{example}

%

\renewcommand\refname{References}
\bibliographystyle{acm}      
\bibliography{YQreference}

\end{document}